\documentclass[10pt]{alggeom}
\usepackage{amsfonts}
\usepackage{mathrsfs}
\usepackage{amscd}
\usepackage{amsmath}
\usepackage{amssymb}
\usepackage{latexsym}
\usepackage{lscape}
\usepackage{comment}
\usepackage{amscd}
\usepackage{wasysym}  
\usepackage{tikz}
\usetikzlibrary{matrix,arrows}
\usepackage{tikz-cd}
\usepackage{appendix}
\usepackage[nice]{nicefrac}
\usepackage{dsfont, stmaryrd}
\usepackage{xparse}

\usepackage[pagebackref,hyperindex,linktocpage=true]{hyperref}
\hypersetup{
    colorlinks,
    linkcolor={red!50!black},
    citecolor={blue!50!black},
    urlcolor={blue!80!black}
}

\usepackage{color}

\usetikzlibrary{decorations.markings}

\makeatletter
\tikzcdset{
  open/.code     = {\tikzcdset{hook, circled};},
  closed/.code   = {\tikzcdset{hook, slashed};},
  open'/.code    = {\tikzcdset{hook', circled};},
  closed'/.code  = {\tikzcdset{hook', slashed};},
  circled/.code  = {\tikzcdset{markwith = {\draw (0,0) circle (.375ex);}};},
  slashed/.code  = {\tikzcdset{markwith = {\draw[-] (-.4ex,-.4ex) -- (.4ex,.4ex);}};},
  markwith/.code ={
    \pgfutil@ifundefined%
    {tikz@library@decorations.markings@loaded}%
    {\pgfutil@packageerror{tikz-cd}{You need to say %
      \string\usetikzlibrary{decorations.markings} to use arrows with markings}{}}{}%
    \pgfkeysalso{/tikz/postaction = {
      /tikz/decorate,
      /tikz/decoration={markings, mark = at position 0.5 with {#1}}}
    }
  },
}
\makeatother

\renewcommand{\geq}{\geqslant}
\renewcommand{\leq}{\leqslant}

\newtheorem{thm}{Theorem}[section]

\newtheorem{propo}[thm]{Proposition}

\newtheorem{propodef}[thm]{Definition and Proposition}
\newtheorem{lem}[thm]{Lemma}
\newtheorem{sublem}[thm]{Sublemma}
\newtheorem{lem-def}[thm]{Lemma-Definition}
\newtheorem{cor}[thm]{Corollary}
\newtheorem{conject}[thm]{Conjecture}
\newtheorem{propert}[thm]{Properties}
\newtheorem{observ}[thm]{Observation}

\newtheorem{fac}[thm]{Fact}

\newtheorem{notat}[thm]{Notation}
\theoremstyle{definition}
\newtheorem*{ack}{Acknowledgement}

\newtheorem{ex}[thm]{Example}

\newtheorem{rmk}[thm]{Remark}
\newtheorem{dfn}[thm]{Definition}
\newtheorem{quest}[thm]{Question}

\newtheorem*{abs}{Abstract}
\numberwithin{equation}{section}

\newcommand{\nc}{\newcommand}

\nc{\abst}{\begin{abs}} \nc{\xabst}{\end{abs}}
\nc{\theo}{\begin{thm}} \nc{\xtheo}{\end{thm}}
\nc{\prop}{\begin{propo}} \nc{\xprop}{\end{propo}}

\nc{\nota}{\begin{notat}} \nc{\xnota}{\end{notat}}
\nc{\depr}{\begin{propodef}} \nc{\xdepr}{\end{propodef}}
\nc{\lemm}{\begin{lem}} \nc{\xlemm}{\end{lem}}
\nc{\sublemm}{\begin{sublem}} \nc{\xsublemm}{\end{sublem}}
\nc{\lemmdefi}{\begin{lem-def}} \nc{\xlemmdefi}{\end{lem-def}}
\nc{\coro}{\begin{cor}} \nc{\xcoro}{\end{cor}}
\nc{\conj}{\begin{conject}} \nc{\xconj}{\end{conject}}
\nc{\proper}{\begin{propert}} \nc{\xproper}{\end{propert}}
\nc{\obse}{\begin{observ}} \nc{\xobse}{\end{observ}}
\nc{\ques}{\begin{quest}} \nc{\xques}{\end{quest}}

\nc{\fact}{\begin{fac}} \nc{\xfact}{\end{fac}}
\nc{\ackn}{\begin{ack}} \nc{\xackn}{\end{ack}}
\nc{\exam}{\begin{ex}} \nc{\xexam}{\end{ex}}
\nc{\rema}{\begin{rmk}} \nc{\xrema}{\end{rmk}}
\nc{\defi}{\begin{dfn}} \nc{\xdefi}{\end{dfn}}

\nc{\pf}{\begin{proof}} \nc{\xpf}{\end{proof}}

\nc{\nn}{\mathrm{n}}
\nc{\on}{\operatorname}
\nc{\fraka}{{\mathfrak a}} \nc{\bba}{{\mathbf a}}
\nc{\frakb}{{\mathfrak b}}
\nc{\frakc}{{\mathfrak c}}
\nc{\frakd}{{\mathfrak d}}
\nc{\frake}{{\mathfrak e}}
\nc{\frakf}{{\mathfrak f}}
\nc{\frakg}{{\mathfrak g}}
\nc{\frakh}{{\mathfrak h}}
\nc{\fraki}{{\mathfrak i}}
\nc{\frakj}{{\mathfrak j}}
\nc{\frakk}{{\mathfrak k}}
\nc{\frakl}{{\mathfrak l}}
\nc{\frakm}{{\mathfrak m}}
\nc{\frakn}{{\mathfrak n}}
\nc{\frako}{{\mathfrak o}}
\nc{\frakp}{{\mathfrak p}}
\nc{\frakq}{{\mathfrak q}}
\nc{\frakr}{{\mathfrak r}}
\nc{\fraks}{{\mathfrak s}}
\nc{\frakt}{{\mathfrak t}}
\nc{\fraku}{{\mathfrak u}}
\nc{\frakv}{{\mathfrak v}}
\nc{\frakw}{{\mathfrak w}}
\nc{\frakx}{{\mathfrak x}}
\nc{\fraky}{{\mathfrak y}}
\nc{\frakz}{{\mathfrak z}}
\nc{\frakA}{{\mathfrak A}}
\nc{\frakB}{{\mathfrak B}}
\nc{\frakC}{{\mathfrak C}}
\nc{\frakD}{{\mathfrak D}}
\nc{\frakE}{{\mathfrak E}}
\nc{\frakF}{{\mathfrak F}}
\nc{\frakG}{{\mathfrak G}}
\nc{\frakH}{{\mathfrak H}}
\nc{\frakI}{{\mathfrak I}}
\nc{\frakJ}{{\mathfrak J}}
\nc{\frakK}{{\mathfrak K}}
\nc{\frakL}{{\mathfrak L}}
\nc{\frakM}{{\mathfrak M}}
\nc{\frakN}{{\mathfrak N}}
\nc{\frakO}{{\mathfrak O}}
\nc{\frakP}{{\mathfrak P}}
\nc{\frakQ}{{\mathfrak Q}}
\nc{\frakR}{{\mathfrak R}}
\nc{\frakS}{{\mathfrak S}}
\nc{\frakT}{{\mathfrak T}}
\nc{\frakU}{{\mathfrak U}}
\nc{\frakV}{{\mathfrak V}}
\nc{\frakW}{{\mathfrak W}}
\nc{\frakX}{{\mathfrak X}}
\nc{\frakY}{{\mathfrak Y}}
\nc{\frakZ}{{\mathfrak Z}}
\nc{\bbA}{{\mathbb A}}
\nc{\bbB}{{\mathbb B}}
\nc{\bbC}{{\mathbb C}}
\nc{\bbD}{{\mathbb D}}
\nc{\bbE}{{\mathbb E}}
\nc{\bbF}{{\mathbb F}} \nc{\bbf}{{\mathbf f}}
\nc{\bbG}{{\mathbb G}}
\nc{\bbH}{{\mathbb H}}
\nc{\bbI}{{\mathbb I}}
\nc{\bbJ}{{\mathbb J}}
\nc{\bbK}{{\mathbb K}}
\nc{\bbL}{{\mathbb L}}
\nc{\bbM}{{\mathbb M}}
\nc{\bbN}{{\mathbb N}}
\nc{\bbO}{{\mathbb O}}
\nc{\bbP}{{\mathbb P}}
\nc{\bbQ}{{\mathbb Q}}
\nc{\bbR}{{\mathbb R}}
\nc{\bbS}{{\mathbb S}}
\nc{\bbT}{{\mathbb T}}
\nc{\bbU}{{\mathbb U}}
\nc{\bbV}{{\mathbb V}}
\nc{\bbW}{{\mathbb W}}
\nc{\bbX}{{\mathbb X}}
\nc{\bbY}{{\mathbb Y}}
\nc{\bbZ}{{\mathbb Z}}
\nc{\calA}{{\mathcal A}}
\nc{\calB}{{\mathcal B}}
\nc{\calC}{{\mathcal C}}
\nc{\calD}{{\mathcal D}}
\nc{\calE}{{\mathcal E}}
\nc{\calF}{{\mathcal F}}
\nc{\calG}{{\mathcal G}}
\nc{\calH}{{\mathcal H}}
\nc{\calI}{{\mathcal I}}
\nc{\calJ}{{\mathcal J}}
\nc{\calK}{{\mathcal K}}
\nc{\calL}{{\mathcal L}}
\nc{\calM}{{\mathcal M}}
\nc{\calN}{{\mathcal N}}
\nc{\calO}{{\mathcal O}}
\nc{\calP}{{\mathcal P}}
\nc{\calQ}{{\mathcal Q}}
\nc{\calR}{{\mathcal R}}
\nc{\calS}{{\mathcal S}}
\nc{\calT}{{\mathcal T}}
\nc{\calU}{{\mathcal U}}
\nc{\calV}{{\mathcal V}}
\nc{\calW}{{\mathcal W}}
\nc{\calX}{{\mathcal X}}
\nc{\calY}{{\mathcal Y}}
\nc{\calZ}{{\mathcal Z}}

\nc{\Spac}{{\mathrm{Spaces}}}
\nc{\scrA}{{\mathscr A}}
\nc{\scrE}{{\mathscr E}}
\nc{\scrR}{{\mathscr R}}

\nc{\Bmu}{\mbox{$\raisebox{-0.59ex}{$l$}\hspace{-0.18em}\mu\hspace{-0.88em}\raisebox{-0.98ex}{\scalebox{2}{$\color{white}.$}}\hspace{-0.416em}\raisebox{+0.88ex}{$\color{white}.$}\hspace{0.46em}$}{}}

\nc{\bnu}{{\bar{ \nu}}}

\nc{\olO}{\bar{\calO}}

\nc{\al}{{\alpha}} 
\nc{\be}{{\beta}}
\nc{\ga}{{\gamma}} \nc{\Ga}{{\Gamma}}
 \nc{\hGa}{\hat{\Gamma}}
\nc{\ve}{{\varepsilon}} 
\nc{\la}{{\lambda}} \nc{\La}{{\Lambda}}
\nc{\om}{\omega} \nc{\Om}{\Omega} 
\nc{\sig}{{\sigma}} \nc{\Sig}{{\Sigma}}

\nc{\tnb}{\psi_{\rm tame}}
\nc{\oM}{\overline{{M}}}
\nc{\op}{{\on{op}}}
\nc{\ad}{{\on{ad}}}
\nc{\alg}{{\on{alg}}}
\nc{\Ad}{{\on{Ad}}}
\nc{\Adm}{{\on{Adm}}} \nc{\aff}{{\on{aff}}}
\nc{\Aut}{{\on{Aut}}}
\nc{\Bun}{{\on{Bun}}}
\nc{\cha}{{\on{char}}}
\nc{\der}{{\on{der}}}
\nc{\Der}{{\on{Der}}}
\nc{\diag}{{\on{diag}}}
\nc{\End}{{\on{End}}}
\nc{\Fl}{{\calF\!\ell}}
\nc{\Tr}{{\on{Transp}}}
\nc{\TR}{{\calT\!\calR}}
\nc{\Gal}{{\on{Gal}}}
\nc{\Gr}{{\on{Gr}}}
\nc{\rH}{{\on{H}}}
\nc{\Hom}{{\on{Hom}}}
\nc{\IC}{{\on{IC}}}
\nc{\id}{{\on{id}}}
\nc{\Id}{{\on{Id}}}
\nc{\ind}{{\on{ind}}}
\nc{\Ind}{{\on{Ind}}}
\nc{\Lie}{{\on{Lie}}}
\nc{\Pic}{{\on{Pic}}}
\nc{\pr}{{\on{pr}}}
\nc{\Res}{{\on{Res}}}
\nc{\res}{{\on{res}}} \nc{\Sat}{{\on{Sat}}}
\nc{\s}{{\on{sc}}}
\nc{\drv}{{\on{der}}}
\nc{\sgn}{{\on{sgn}}}
\nc{\Spec}{{\on{Spec}}}\nc{\Spf}{\on{Spf}} 
\nc{\Sph}{\on{Sph}}
\nc{\St}{{\on{St}}}
\nc{\tr}{{\on{tr}}}
\nc{\Mod}{{\mathrm{-Mod}}}
\nc{\Hilb}{{\on{Hilb}}} 
\nc{\Ext}{{\on{Ext}}} 
\nc{\vs}{{\on{Vec}}}
\nc{\ev}{{\on{ev}}}
\nc{\nO}{{\breve{\calO}}}
\nc{\tS}{{\tilde{S}}}
\nc{\spe}{{\on{sp}}}
\nc{\loc}{{\on{loc}}}
\nc{\Sym}{{\on{Sym}}}
\nc{\Cone}{{\on{C}}}
\nc{\syn}{{\on{syn}}}
\nc{\reg}{{\on{reg}}}
\nc{\colim}{{\on{colim}}}
\nc{\Norm}{{\on{N}}}

\nc{\nscrR}{{\mathscr{R}^{\on{nr}}}}

\nc{\GL}{{\on{GL}}}
\nc{\U}{{\on{U}}}
\nc{\Gl}{\on{Gl}} 
\nc{\GSp}{{\on{GSp}}}
\nc{\gl}{{\frakg\frakl}}
\nc{\SL}{{\on{SL}}} 
\nc{\SU}{{\on{SU}}} 
\nc{\SO}{{\on{SO}}}
\nc{\PGL}{{\on{PGL}}}

\nc{\Conv}{{\on{Conv}}}
\nc{\Rep}{{\on{Rep}}}
\nc{\Dom}{{\on{Dom}}}
\nc{\red}{{\on{red}}}
\nc{\act}{{\on{act}}}
\nc{\nr}{{\on{nr}}}
\nc{\ctf}{{\on{ctf}}}

\nc{\str}{{\on{-}}} 
\nc{\os}{{\bar{s}}}
\nc{\oeta}{{\bar{\eta}}}

\nc{\hookto}{\hookrightarrow}
\nc{\longto}{\longrightarrow}
\nc{\leftto}{\leftarrow}
\nc{\onto}{\twoheadrightarrow}
\nc{\lonto}{\twoheadleftarrow}
\newcommand*\isomto{%
  \renewcommand{\arraystretch}{0.1}
  \begin{array}[b]{c} {}_{\sim} \\ \longrightarrow \end{array}%
}

\nc{\uG}{{\underline{G}}}
\nc{\uA}{{\underline{A}}}
\nc{\uS}{{\underline{S}}}
\nc{\uT}{{\underline{T}}}
\nc{\uM}{{\underline{M}}}
\nc{\uP}{{\underline{P}}}
\nc{\uB}{{\underline{B}}}
\nc{\uN}{{\underline{N}}}

\nc{\ucG}{{\underline{\calG}}}
\nc{\ucA}{{\underline{\calA}}}
\nc{\ucS}{{\underline{\calS}}}
\nc{\ucT}{{\underline{\calT}}}
\nc{\ucalM}{{\underline{\calM}}}
\nc{\ucP}{{\underline{\calP}}}
\nc{\ucalN}{{\underline{\calN}}}

\nc{\bF}{{\breve{F}}}

\nc{\oFl}{{\overline{\Fl}}} 
\nc{\bU}{{\overline{U}}}
\nc{\tGr}{{\tilde{\Gr}}}
\nc{\cGr}{\calG\! r}
\nc{\oGr}{\overline{\on{Gr}}} 
\nc{\ocGr}{\overline{\calG\! r}}
\nc{\co}{{\colon}}
\nc{\sch}[1]{(Sch/{#1})}
\nc{\HypLoc}[1]{HypLoc({#1})}

\nc{\ohtimes}{\stackrel{!}{\otimes}}
\nc{\boxtilde}{\widetilde{\boxtimes}}
\nc{\vstar}{{\varhexstar}}

\nc{\Div}{\on{Div}}
\nc{\Sht}{\on{Sht}}
\nc{\Frob}{\on{Frob}}

\nc{\x}{\times}
\nc{\bsl}{\backslash}
\nc{\algQl}{{\bar{\bbQ}_\ell}}
\nc{\sF}{{\bar{F}}}
\nc{\nF}{{\breve{F}}}
\nc{\nW}{{W^{\on{nr}}}}
\nc{\sk}{{\bar{k}}}
\nc{\cont}{\on{c}}
\nc{\Supp}{\on{Supp}}
\nc{\blt}{\bullet}  
\nc{\dom}{\on{dom}}
\nc{\scon}{{\on{sc}}} 
\nc{\Affine}{\on{Aff}} 
\nc{\nscrA}{\mathscr{A}^{\on{nr}}} 
\nc{\nfraka}{{\bbf^{\on{nr}}}}
\nc{\ran}{{\rangle}}
\nc{\lan}{{\langle}}
\nc{\bk}{{\bar{k}}}
\nc{\tF}{{\tilde{F}}}
\nc{\sS}{{\bar{S}}}
\nc{\LG}{{^\text{L}\hspace{-0.04cm}G}}
\nc{\LL}{{^\text{L}\hspace{-0.07cm}L}}
\nc{\et}{{\text{\rm \'et}}}
\nc{\inv}{{\on{inv}}}
\nc{\Hecke}{{\on{Hecke}}}
\nc{\Isom}{{\on{Isom}}}
\nc{\oSht}{{\overline{\on{Sht}}}}
\nc{\umu}{{\underline \mu}}
\nc{\AIJ}{{\calO_X[{\scriptstyle{\calI\over \calJ}}]}}
\nc{\Proj}{{\on{Proj}}}
\nc{\Bl}{{\on{Bl}}}

\nc{\Pos}{{\on{Pos}}}
\nc{\Sets}{{\on{Sets}}}
\nc{\AffSch}{{\on{AffSch}}}
\nc{\Groups}{{\on{Groups}}}
\nc{\Gpds}{{\on{Groupoids}}}
\nc{\Sch}{{\on{Sch}}}
\nc{\fl}{{\on{flat}}}

\nc{\pot}[1]{ [\hspace{-0,5mm}[ {#1} ]\hspace{-0,5mm}] }
\nc{\rpot}[1]{ (\hspace{-0,7mm}( {#1} )\hspace{-0,7mm}) }

\nc{\defined}{\hspace{0.1cm}\stackrel{\text{\tiny \rm def}}{=}\hspace{0.1cm}}

\begin{document}

\pagestyle{plain}

\title{Dilatations of categories}

\classification{18E35}
\keywords{dilatations of categories, categories of fractions, localizations of categories, algebraic dilatations}

\maketitle
\[\text{{\Large Arnaud Mayeux     }}\]

\[\text{\today} \]

$~~$

\abst
Dilatations modify categories by imposing that some morphisms factorize through some others. This is formalized by a universal property. This text is devoted to introduce and study this construction.  Examples of dilatations of categories include localizations of categories and dilatations of rings.
\xabst

\tableofcontents

\ackn
 This project has received funding from the ISF grant 1577/23 and the ERC grant agreement No 101002592. I am grateful to the referee for his suggestion to remove an unnecessary assumption and other comments. 
\xackn

\section{Introduction} Given a subset $S$ of a commutative ring $A$, one has the localization $S^{-1}A$ of $A$ relatively to $S$. Needless to insist on the fact that this is a fundamental construction. The localization process of a commutative ring is extended in several ways, among them:
\begin{enumerate}
\item localization of categories (cf. e.g. \cite{GM}, \cite{GZ} and \cite{B}), a basic construction used in many branchs of mathematics;
\item dilatation of rings (cf. e.g. \cite[§2]{M} and \cite[\href{https://stacks.math.columbia.edu/tag/052P}{Tag 052P}]{stacks-project}), the building blocks of dilatations of schemes (cf. e.g. \cite{DMdS} and the large number of references therein).
\end{enumerate}
In this text, we provide a construction unifying all of these constructions in a single construction: dilatations of categories. Intuitively, localization is a process that adds and imposes inverses of elements or morphisms. In other words, localization adds some "fractions" with prescribed denominators. Intuitively, dilatations do the same thing except that it only adds some fractions where both numerators and denominators are prescribed. 
Let $\calC$ be a category, and let $\Sigma $ be a collection of morphisms of $\calC$. Let us recall some properties of the localization $\calC[\Sigma ^{-1}]$. We have a functor $L:\calC \to \calC[\Sigma^{-1}]$. The objects of the localization $\calC [\Sigma ^{-1}]$ coincide with the objects of $\calC$, i.e. $L$ is the identity on objects. Given a morphism $d$ in $\Sigma $, $L(d)$ is an isomorphism in $\calC[\Sigma ^{-1}]$. If $F: \calC \to \calD$ is another functor such that $F(d) $ is an isomorphism for all $d$ in $\Sigma$, then $F$ factors through $L$.
 Now assume that for any $d$ in $\Sigma$ we have a sieve $N_d$, in $\calC$ over the codomain of $d$. The dilatation process will provide a category $\calC '$ and a functor $\Theta : \calC \to \calC'$. The objects of $ \calC'$ will also coincide with the objects of $\calC$. For any $d \in \Sigma$ and any $n \in N_d$, there exists a unique arrow $b$ in $\calC '$ such that the diagram 
\[ \begin{tikzcd}
X \ar[rr, "{\Theta (n)}"] \ar[rd, "{b}", dotted] & & Y \\
& Z \ar[ru, "{\Theta (d)}" ] &
\end{tikzcd}\] commutes. The element $b$ is thought of as a non-commutative fraction $d \backslash n = d^{-1} \circ n $. Dilatations also satisfy a natural universal property (Theorem \ref{uni}), namely $\Theta $ is universal among all functors $F : \calC \to \calD$ such that 
\begin{enumerate}
\item the sieve generated by $F(N_d) $ is included in the sieve generated by $F(d)$, for all morphisms $d$ in the collection $\Sigma$,
\item the localization $\calD \to \calD [\Sigma ^{-1}]$ is faithful. 
\end{enumerate}
Condition (i) allows existence of the factorization morphisms "b" as before and condition (ii) allows unicity. 
To prove the existence of dilatations in this note, we use directed graphs.

We now discuss the structure of this paper.
Definition \ref{depr1} and antecedent material provide the definition of dilatations of categories. It holds in a very general framework. Proposition \ref{propdil}, Theorem \ref{uni} and Proposition \ref{repre} are the main results. Facts \ref{locasdil} and \ref{Ma} make connections with localizations of categories and dilatations of rings.  We also discuss a related notion in §\ref{codi}, that of codilatations of categories. After that, using the fact that algebraic structures (sets with composition laws and axioms) are in particular categories with a single objects, we introduce dilatations of some non-commutative structures (e.g. monoids).
Finally, Fact \ref{contrexemplecat} shows that the "natural generalization" of another caracterization of dilatations of rings, in terms of localizations, does not hold for categories.

\section{Dilatations of categories: definition via fractions}

We now fix a category  $\calC$ with objects $Ob\calC $ and morphisms $Mor \calC$.

\subsection{Sieves}

 Recall that if $X\xrightarrow{a}Y$ belongs to $Mor \calC$, $dom(a):= X$ is called the domain and $cod(a):=Y$ the codomain (note that $dom$ and $cod$ also make sense for directed graphs). Recall that a sieve in $\calC$ is a collection of morphisms stable by precompositions. Note that a union of sieves is a sieve.  A sieve over an object $Y$ is a sieve made of morphisms  with codomain $Y$. Similarly, a cosieve is a collection of morphisms stable by postcompositions. 

\defi Let $E$ be a collection of morphisms of $\calC$. The collection $e \circ  n $ with $n \in Mor \calC $, $e \in E$ and $cod(n)=dom(e)$ is a sieve called the sieve generated by $E$. 
This sieve is denoted $S^{\calC}_{E}$. If $E= \{d\}$ is a singleton, we put $S^{\calC}_{E}=S^{\calC}_{d}$. Note that if $d = Id _X$, $S^{\calC } _{Id_{X}}$ is the sieve made of all morphisms with codomain $X$. Similarly, we use the notation $CoS^{\calC}_{E}$ to denote the cosieve generated by $E$.
\xdefi

\subsection{Directed graphs and localizations}
\label{subsectgraph}
 In this section, we discuss localizations of categories.  The content of this subsection is classical and does not claim originality (cf. e.g. \cite{GZ}, \cite{GM} and \cite{B}), however we provide a self-contained description of this construction. Let $\Sigma $ be a collection of morphisms of $\calC$.

\defi \label{graphdef} 
Let $\calG $ be the oriented graph defined as follows. The vertices of $\calG$ are equal to the objects of $\calC$. The directed lines of $\calG$ are made of \begin{enumerate} \item for each morphism $a$ of $\calC$, a directed line $dom(a) \xrightarrow{a} cod(a)$ of $\calG $,
\item for each morphism $d$ in $\Sigma $ a directed line $cod(d) \xrightarrow{l_d} dom(d)$ of $\calG$ (in particular $dom (l_d) = cod(d) $ and $cod (l_d) = dom (d)$).
\end{enumerate}
\xdefi

\defi
A $\Sigma$-sequence of directed lines in $\calG$ is a sequence, finite and possibly empty, of directed lines $x_1, x_2 , \ldots , x_n$ of $\calG$ such that  $cod(x_i)= dom (x_{i+1})$ ($n\geq0$ is an integer). By convention, an empty sequence is just an element in the collection $Ob \calC$. For a non-empty sequence $s=(x_1 , \ldots , x_n )$, we define the domain as $dom(s )= dom (x_1) $ and the codomain as $cod (s) = cod (x_n)$. If a sequence $s$ is empty and given by an object $X$, then we put $dom(s)=cod(s)=X$. Note that $\Sigma$-sequences of directed lines with compatible domains and codomains can be composed (by convention, composing a sequence $s$ with an empty sequence $e$ gives $s$).
\xdefi

\defi  \label{sequencedefilict}
We say that two  $\Sigma$-sequences of directed lines $s$ and $s'$ in $\calG$ are equivalent if $dom(s)=dom(s'),$ $ cod(s)=cod(s')$ and if one can be obtained from the other by a chain of elementary equivalences of the following types:
\begin{enumerate}
\item a sequence  $x_1, \ldots , x_n $ such that  $x_i , x_{i+1}$  are equal to $a,a' \in Mor \calC$ and $cod(a)=dom(a')$, for some $1 \leq i , i+1 \leq n$, is equivalent to the sequence $x_1 , \ldots , x_{i-1} , x' , x_{i+2}, \ldots , x_{n}$ where $x' $ is the composition of $a $ and $a'$ in $\calC$,
\item a sequence $x_1 , \ldots , x_n $ such that $x_i , x_{i+1} $ are equal to $ l_d , d$ with $d \in \Sigma$, for some $1 \leq i , i+1 \leq n$, is equivalent to the sequence $x_1, \ldots , x_{i-1} , x_{i+2}, \ldots , x_n $,
\item a sequence $x_1 , \ldots , x_n$ such that $x_i , x_{i+1}$ are equal to $d,l_d$ with $d \in \Sigma $, for some $1 \leq i , i+1 \leq n$, is equivalent to the sequence $x_1, \ldots , x_{i-1} , x_{i+2}, \ldots , x_n$,
\item for any object $X$, the sequence $Id_X$ is equivalent to the empty sequence at $X$.
\end{enumerate} 
In other words and informally, two $\Sigma$-sequences of directed lines are equivalent if one can be obtained from the other  by the operations consisting in exchanging parts of sequences as follows: \begin{align*}
  X \xrightarrow{a} Y \xrightarrow{a'} Z  & \leftrightsquigarrow  X \xrightarrow{a' \circ a } Z \\    Y \xrightarrow{l_d } Z \xrightarrow{d} Y & \leftrightsquigarrow  Y  \text{ (the empty sequence at $Y$)} \\   Z \xrightarrow{d} Y \xrightarrow{l_d} Z & \leftrightsquigarrow   Z  \text{ (the empty sequence at $Z$)} \\ X \xrightarrow{Id_X} X & \leftrightsquigarrow X  \text{ (the empty sequence at $X$)}.
\end{align*} The equivalence class of a sequence $s$ is denoted by $[s]$.
Note that equivalence classes of $\Sigma$-sequences of directed lines with compatible domains and codomains can be composed associatively. 
\xdefi

\defi A $\Sigma$-fraction is an equivalence class of $\Sigma$-sequence.

\xdefi

\fact Let $d,d' $ in $\Sigma $ such that $cod(d)=dom(d')$ and let $d''$ be their composite. Assume that $d''$ belongs to $\Sigma$. Then \[[cod(d') \xrightarrow{l_{d'}} dom(d') \xrightarrow{l_{d}} dom (d) ] = [cod (d') \xrightarrow{l_{d''}} dom(d) ].\]
\xfact

\pf This follows from the equalities \begin{align*} &=[cod(d') \xrightarrow{l_{d'}} dom(d') \xrightarrow{l_{d}} dom (d) ]\\ &=[cod(d') \xrightarrow{l_{d'}} dom(d') \xrightarrow{l_{d}} dom (d) \xrightarrow{d''} cod(d') \xrightarrow{l_{d''}} dom (d) ] \\
&= [cod(d') \xrightarrow{l_{d'}} dom(d') \xrightarrow{l_{d}} dom (d) \xrightarrow{d} cod(d)\xrightarrow{d'} cod(d') \xrightarrow{l_{d''}} dom (d) ]  \\ &= [cod (d') \xrightarrow{l_{d''}} dom(d) ] .\end{align*}
\xpf 

\defi The localization of $\calC $ relatively to $\Sigma$ is the category $\calC [\Sigma^{-1}]$ whose objects are the objects of $\calC$ and whose morphisms are $\Sigma$-fractions. We have a canonical functor $L:\calC \to \calC[\Sigma^{-1}]$.
\xdefi

\prop[(Universal property of localization)]
Let $F: \calC \to \calD$ be a functor such that $F(d)$ is invertible for all $d $ in $\Sigma$. Then there exists a unique functor $F': \calC[\Sigma ^{-1}] \to \calD$ such that  $F= F' \circ L$.
\xprop
\pf
We first prove unicity. Assume that $F'$ exists. Let $X$ be in $Ob\calC = Ob \calC '$, then $F'(X) = F'(L(X))=F(X)$. Let $a$ be in $Mor \calC$, then $F'([a])=F'(L(a)) = F(a)$. Let $d $ be in $\Sigma$, then $F'([l_d]) \circ F(d) =F'( [l_d] \circ [d])= Id_{dom(F(d))}$ and $F(d) \circ F'([l_d]) = F' ( [d ]\circ [ l_d  ] ) =Id_{cod(F(d))}$. This implies that $F'([l_d])= F(d)^{-1}$. This proves unicity. To prove existence, it is enough to prove that the assignements $X \mapsto X, [a] \mapsto F(a), l_d \mapsto F(d)^{-1}$ (as before) provide a well-defined functor $F'$. For this it is enough to prove that the assignements are compatible with the elementary equivalences of Definition \ref{sequencedefilict}, which is immediate.
\xpf

\subsection{Dilatations of categories: definition}
\label{subsectdef}

Let $\calC$ be a category. 

\defi
A center in $\calC$ is a collection  $\{[N_i, d_i ]\}_{i \in I }$ of pairs $[N_i , d_i]$, indexed by a collection $I$ and such that, for all $i$ in $I$, $d_i$ is a morphism of $
\calC$ and $N_i$ is a sieve over $cod(d_i)$.
\xdefi

We now fix a center $\{[N_i, d_i ]\}_{i \in I }$ in $\calC$ (we sometimes use the notation $N_{d_i}$ to denote $N_i$). 

Put $\Sigma =\{d_i\}_{i \in I}$.

\defi \label{fractiondef} A $\{[N_i, d_i ]\}_{i \in I }$-fraction is a $\Sigma$-fraction such that a representative can be written as
\begin{center}$ X_1 \xrightarrow{n_1} Y_1 \xrightarrow{l_{d_{i_1}}} X_2 \xrightarrow{n_2} Y_2 \xrightarrow{l_{d_{i_2}}} X_3 \ldots X_{k} \xrightarrow{n_k} Y_k \xrightarrow{l_{d_{i_k}}} X_{k+1} \xrightarrow{a} X_{k+2}$\end{center}  with $a \in Mor \calC$, $k \geq 0$ and ${i_j} \in I$, $n_i \in N_{i_j} $ for all $j \in \{1 , \ldots , k \}.$
\xdefi

\fact  \label{compsotion}
 $\{[N_i, d_i ]\}_{i \in I }$-fractions with compatible domains and codomains can be composed associatively.
\xfact
\pf
This is immediate since the $N_i's$ are sieves. Associativity is immediate. 
\xpf 

\rema \begin{enumerate} \item If $X\xrightarrow{a}Y$ is a morphism in $\calC$, the $\{d_i\}_{i \in I}$-fraction $[ X\xrightarrow{a}Y ]$ is a $\{[N_i, d_i ]\}_{i \in I }$-fraction. In particular, the class $[Id_X] $ of the identity at an object $X$ is a $\{[N_i, d_i ]\}_{i \in I }$-fraction. \item The $\Sigma$-fraction \[ X_1 \xrightarrow{n_1} Y_1 \xrightarrow{l_{d_{i_1}}} X_2 \xrightarrow{n_2} Y_2 \xrightarrow{l_{d_{i_2}}} X_3 \ldots X_{k} \xrightarrow{n_k} Y_k \xrightarrow{l_{d_{i_k}}} X_{k+1}\] is a $\{[N_i, d_i ]\}_{i \in I }$-fraction since it is equivalent to 
\[ X_1 \xrightarrow{n_1} Y_1 \xrightarrow{l_{d_{i_1}}} X_2 \xrightarrow{n_2} Y_2 \xrightarrow{l_{d_{i_2}}} X_3 \ldots X_{k} \xrightarrow{n_k} Y_k \xrightarrow{l_{d_{i_k}}} X_{k+1} \xrightarrow{Id_{X_{k+1}}} X_{k+1}.\] \end{enumerate} 
\xrema

\defi \label{depr1} The dilatation of $\calC$ with center $\{[N_i, d_i ]\}_{i \in I }$ is the category $\calC'$ defined as follows. The objects of $\calC'$ are equal to the objects of $\calC$. If $X $ and $Y$ are objects in $\calC '$, morphisms between $X $ and $Y$ are given by $\{[N_i, d_i ]\}_{i \in I }$-fractions (cf. Definition \ref{fractiondef}) with domain $X$ and codomain $Y$. Fact \ref{compsotion} shows that $\calC'$ is indeed a category, which is also denoted $\calC \left[ \left\{ (d_i)^{-1} \circ  N_i \right\}_{i \in I} \right]$.
\xdefi

\fact \label{dildansloc}
We have a faithful functor $\calC \left[ \left\{ (d_i)^{-1} \circ N_i \right\}_{i \in I} \right] \to \calC \left[ \left(\left\{ d_i \right\}_{i \in I} \right) ^{-1} \right]$.
\xfact 
\pf
The functor is obtained regarding $\{[N_i, d_i ]\}_{i \in I }$-fractions as $\Sigma$-fractions. 
\xpf

\fact \label{locasdil} Assume that, for all $i $ in $I$, $N_i $ is the sieve  $S_{Id_{cod(d_i)}}^{\calC}$, then $\calC \left[ \left\{ (d_i)^{-1} \circ N_i \right\}_{i \in I} \right]$ is the localization  $\calC [(\{d_i\}_{i \in I}) ^{-1}]$.
\xfact 
\pf In this case, the faithful functor of Fact \ref{dildansloc} is also full, as any $\{d_i\}_{i \in I}$-fraction is equivalent to a $\{[S_{Id_{cod(d_i)}}^{\calC}, d_i ]\}_{i \in I }$-fraction. \xpf 

\coro Assume that $\calC $ is small, then $\calC \left[ \left\{ (d_i)^{-1} \circ N_i \right\}_{i \in I} \right]$ is small.
\xcoro
\pf
Combine Fact \ref{dildansloc} and \cite[Proposition 5.2.2]{B}.
\xpf

\section{Universal property of dilatations and related results}

We proceed with the notation from §\ref{subsectdef}. Recall that $\calC' =  \calC \left[ \left\{ (d_i)^{-1} \circ N_i \right\}_{i \in I} \right]$.

\prop \label{propdil} The following assertions hold. \begin{enumerate} \item We have a canonical functor $\Theta : \calC \to \calC'$ defined as follows.
If $X$ is an object of $\calC$, $\Theta (X) =X$. If $X \xrightarrow{a}Y$ is a morphism of $\calC$, $\Theta (a)= [X \xrightarrow{a}Y]$.
\item For any $i$ in $I$ and any $n \in N_i $, there exists a unique morphism $dom(n) \xrightarrow{b} dom(d_i)$ in $\calC'$ such that the following triangle commutes 
\begin{tikzcd}
dom(n) \ar[rr, "{[n]}"] \ar[rd, "{b}"] & & cod(n) \\
& dom(d_i) \ar[ru, "{[d_i]}" ] &
\end{tikzcd}.  
 Intuitively, the morphism $b$ is thought of as a non-commutative fraction $d_i \backslash n = (d_i)^{-1} \circ n $ of morphisms. Mathematically, we have \[b = [ X \xrightarrow{n} Y \xrightarrow{l_{d_i}} Z]. \]
\end{enumerate}
\xprop 

\pf \begin{enumerate} \item Obvious.
\item Put $X= dom(n), Y=cod(n)$ and $Z = dom(d_i)$.
The class $[X \xrightarrow{n} Y \xrightarrow{l_{d_i}} Z]$ composed with the class  $[Z \xrightarrow{d_i} Y]$ equals the class $[n]$. This shows that there exists an arrow $b$ in $\calC'$ such that the triangle commutes. We now prove unicity. Assume that  $m \in Mor \calC'$ satisfies that the composition $ X \xrightarrow{m} Z \xrightarrow{\Theta(d_i)}  Y $ equals $ X \xrightarrow{\Theta (n) } Y$. Let $s$ be a representative sequence of $m$. Then the class  $ [X \xrightarrow{s} Z \xrightarrow{d_i}  Y \xrightarrow{l_{d_i}} Z]$ is equal to the class $ [X \xrightarrow{n  } Y \xrightarrow{l_{d_i}}Z]$ and also to the class $[ X \xrightarrow{s} Z ]$. This shows that $m$ equals  $[X \xrightarrow{n} Y \xrightarrow{l_{d_i}} Z]$. 
\end{enumerate}
\xpf

Recall that a bimorphism is defined as a morphism which is both a monomorphism and an epimorphism. Bimorphisms are also called regular morphisms. 
\fact \label{bimo} A morphism whose image under a faithful functor is a bimorphism is itself a bimorphism.
\xfact 
\pf
Let $F :\calA \to \calB$ be a faithful functor. Let $f$ be a morphism of $\calA$ such that $F(f)$ is a bimorphism. Let $a,b,c,d$ be morphisms of $\calA$ such that $a \circ f = b \circ f$ and $f \circ c = f \circ d$. Then $F(a) \circ F(f) =F( b) \circ F(f)$ and $F(f) \circ F(c) = F(f) \circ F(d)$. So $F(a)=F(b) $ and $F(c)=F(d)$. Consequently $a=b$ and $c=d$. This proves that $f$ is a bimorphism. 
\xpf 
\prop \label{nonz} Let $i \in I$, then $ \Theta(d_i)$ is a bimorphism in $\calC'$.
\xprop
\pf A morphism whose image under a faithful functor is a bimorphism is itself a bimorphism by Fact \ref{bimo}. Now Proposition \ref{nonz} follows from Fact \ref{dildansloc}.
\xpf

\defi For any $i \in I $, let $S_{\Theta (N_i)}^{ \calC ' }$ be the sieve over $cod(\Theta (d_i))$ in $\calC'$ generated by $\{\Theta (n) | n \in N_i \}.$  Similarly, let $S^{\calC'}_{\Theta (d_i)}$ be the sieve over $cod (\Theta (d_i))$ generated by $\Theta (d_i)$.
\xdefi

\prop \label{exc} Let $i \in I  $,  then  $ S_{\Theta (N_i)}^{ \calC ' } \subset  S^{\calC'}_{\Theta (d_i)} $.
\xprop 
\pf Let $x $ be in $ S_{\Theta (N_i)}^{ \calC ' }  $. Then $x = [n] \circ t $ with $n \in N_i$. We have \[  x  = [n] \circ t  = [d_i] \circ [ dom (n) \xrightarrow{n} cod(n) \xrightarrow {l_{d_i}} dom(d_i) ] \circ t  .\] This finishes the proof.
\xpf

\defi   Let $Cat _{\calC}^{\Sigma \text{-reg}}$ be the full subcategory of
the comma category $\calC/Cat$ whose objects $F: \calC \to \calD$ are arbitrary functors out of $\calC$ such that the localization functor $\calD \to \calD[F(\Sigma)^{-1}]$ is faithful.
\xdefi

\fact \label{obecjtreg}
The functor $ \Theta: \calC \to  \calC \left[ \left\{ (d_i)^{-1} \circ N_i \right\}_{i \in I} \right]$  is an object in the category $Cat _{\calC}^{\Sigma \text{-reg}}$.
\xfact

\pf \begin{sloppypar}Fact \ref{dildansloc} says that $ \calC \left[ \left\{ (d_i)^{-1} \circ N_i \right\}_{i \in I} \right] \to  \calC [( \{d_i\}_{i \in I}) ^{-1}]$ is faithful.
Now we observe that $ \calC \left[ \left\{ (d_i)^{-1} \circ N_i \right\}_{i \in I} \right] \left[ (\{\Theta ( d_i) \}_{i \in I})^{-1} \right] $ identifies with the localization $ \calC [( \{d_i\}_{i \in I}) ^{-1}]$. This finishes the proof.  \end{sloppypar}
\xpf 

\fact \label{factr}
Let $\Sigma' \subset \Sigma $ be a subcollection. If $L: \calC \to \calC [\Sigma ^{-1}]$ is faithful, then $L' : \calC \to \calC [ \Sigma '^{-1} ]$ is also faithful.
\xfact
\pf
 Observe that $\calC [\Sigma ^{-1}]= \calC [\Sigma '^{-1}] [ L' ( \Sigma )^{-1}]$, so that we have a commutative triangle of functors \[\begin{tikzcd} \calC \ar[dr, "L"] \ar[d, "L'"] \\ \calC [\Sigma '^{-1}] \ar[r, "l"]& \calC [\Sigma ^ {-1}] \end{tikzcd}.\] Now let $a,b$ be two morphisms such that $L'(a)=L'(b)$. We obtain $l(L'(a))=l(L'(b))$ and $a=b$. So $L'$ is faithful.
\xpf
\rema
We remark that \cite[Lemma 4.4]{BM} shows that if $\calC$ is semi-abelian and integral (cf. \cite{BM} and references therein for precise definitions) and if $B $ is the collection of all bimorphisms of $\calC$, then the identity functor $\calC \to \calC$ itself belongs to $Cat _{\calC}^{B \text{-reg}}$. Still in the above setting of \cite{BM}, \cite[Lemma 4.4]{BM} together with Fact \ref{factr} implies that the identity functor $\calC \to \calC$  itself belongs to $Cat _{\calC}^{\Sigma \text{-reg}}$ for any subcollection $\Sigma $ of $B$. Note that we do not use \cite[Lemma 4.4]{BM} in the present paper.
\xrema

\theo\label{uni} (Universal property)  Let $F: \calC \to \calD$ be an object in $Cat _{\calC}^{\Sigma \text{-reg}}$ such that
for any $i $ in $I$,  we have \[S^{\calD}_{F(N_i)} \subset S^{\calD}_{F(d_i)} .\]
Then there is a unique functor $F' : \calC ' \to \calD$ such that the triangle of functors \[
\begin{tikzcd}
\calC \ar[rr, "F"] \ar[rd, "\Theta "] & & \calD \\
& \calC ' \ar[ru, " F' " ] &
\end{tikzcd}  \] commutes (recall that $\calC '=\calC \left[ \left\{ (d_i)^{-1} \circ N_i \right\}_{i \in I} \right]$). 
\xtheo 
\pf
Assume that such a $F'$ exists, then  $F'(X) = F' (\Theta (X)) = F(X) $ for all $X$ in $Ob \calC  = Ob \calC'$. Now any morphism in $\calC'$ is a composition of morphisms of the form $[X \xrightarrow{n} Y \xrightarrow{l_{d_i}} Z]$ with $i \in I, n \in N_i$, or of the form $[dom(a) \xrightarrow{a} cod(a)]$. Necessarily we have $F' ( [ X \xrightarrow{n} Y \xrightarrow{l_{d_i}} Z ] ) = t$ where $t$ is the unique morphism such that \begin{tikzcd}F(X) \ar[rr, "F(n)"] \ar[rd, "t"] & & F(Y) \\& F(Z) \ar[ru, " F (d_i) " ] &\end{tikzcd} commutes (this $t$ exists since $F(n) $ belongs to $S_{F(d_i)}^{\calD}$ by assumption and is unique because $\calD \to \calD[F(\Sigma) ^{-1}]$ is faithful and so $F(d_i)$ is a bimorphism by Fact \ref{bimo}). Necessarily we have $F'([a]) = F' ( \Theta (a)) = F(a)$.  This shows that $F'$ is unique and given by the formula above if it exists. Reciprocally, the assignement $X \mapsto F(X) , [a] \mapsto F(a), [ X \xrightarrow{n} Y \xrightarrow{l_{d_i}} Z ] \mapsto t \text{ (the unique $t$ as before)}$, is a well-defined functor $F'$ satisfying that $F= F' \circ \Theta$.
\xpf

\fact  \label{factbd}
Let $H: \calA \to \calB $ be a functor between two categories. Let $ E $ be a collection of morphisms of $\calA$. Then $S^{\calB}_{H(E)} = S^\calB _{ H (S^{\calA} _{E} ) }$.
\xfact
\pf
We have $E \subset S^{\calA} _{E}$, so $S^{\calB}_{H(E)} \subset S^\calB _{ H (S^{\calA} _{E} ) }$. Reciprocally, let $x$ be in $S^\calB _{ H (S^{\calA} _{E} ) }$, then there exists $b \in Mor \calB $,  $a \in Mor \calA$ and $e \in E$, such that $x = H ( e \circ a ) \circ b$. Since $H$ is a functor, we have $x = H(e) \circ  ( H(a) \circ b) $. So $x  $ belongs to $S^{\calB}_{H(E)}$. 
\xpf

\prop  \label{repre} The functor $\Theta : \calC \to  \calC \left[ \left\{ (d_i)^{-1} \circ N_i \right\}_{i \in I} \right]  $ represents the covariant functor 
 $Cat_{\calC}^{\Sigma\text{-reg}}\to Set$
given by
\begin{equation}
(\calC \xrightarrow{F} \calD) \;\longmapsto\; \begin{cases}\{*\}, \; \text{ if $S^{\calD}_{F(N_i)} \subset S^{\calD}_{F(d_i)}$ for any $i $ in $I$;}\\ \varnothing,\;\text{else.}\end{cases}
\end{equation}
\xprop
\pf
Let $(\calC \xrightarrow{F} \calD)$ be an object of $Cat_{\calC}^{\Sigma\text{-reg}}$. If $F' \in Hom _{\calC} ( \calC' , \calD) $, then for all $i$ in $I$ we have \[S_{F(N_i)}^{\calD} = S_{F' (\Theta (N_i))}^{\calD} = S_{F' (S_{\Theta (N_i)}^{\calC'})}^{\calD} \subset S_{F' (S_{\Theta(d_i)}^{\calC'})}^{\calD} = S_{F' (\Theta(d_i))}^{\calD}= S_{F(d_i)}^{\calD}.\]
We used Proposition \ref{exc} and Fact \ref{factbd}.
Therefore, if there exists $i$ such that $S^{\calD}_{F(N_i)} \not \subset S^{\calD}_{F(d_i)}$, then $ Hom _{\calC} ( \calC' , \calD) = \varnothing$. If $S^{\calD}_{F(N_i)} \subset S^{\calD}_{F(d_i)}$ for all $i $ in $I$, Theorem \ref{uni} implies that  $ Hom _{\calC} ( \calC' , \calD) $ is a singleton $\{*\}$. This finishes the proof.
\xpf

\fact \label{inclusionfonctoriel} For each $i$ in $I$, let $M_i$ be a subsieve of $N_i$. Then we have a canonical functor 
\[ \varphi: \calC \left[ \left\{ (d_i)^{-1} \circ M_i \right\}_{i \in I} \right] \to\calC \left[ \left\{ (d_i)^{-1} \circ N_i \right\}_{i \in I} \right] . \] The functor $\varphi$ is faithful.
\xfact
\pf Any $\{[M_i, d_i ]\}_{i \in I }$-fraction is a $\{[N_i, d_i ]\}_{i \in I }$-fraction since $M_i $ is a subsieve of $N_i$ for all $i$.
\xpf 

\prop \label{fi}Let $K \subset I$ be a subcollection. Then we have a canonical functor 
\[ \Phi: \calC \left[ \left\{ (d_i)^{-1} \circ N_i \right\}_{i \in K} \right]  \to \calC \left[ \left\{ (d_i)^{-1} \circ N_i \right\}_{i \in I} \right] . \]
Moreover \begin{enumerate} 
\item if $N_i = S ^{\calC} _{d_i}$ for all $i $ in $I \smallsetminus K $, then $\Phi$ is full,
 \item if $\calC [ (\{d_i\}_{i \in K})^{-1} ] \to  \calC [ (\{d_i\}_{i \in I})^{-1 } ] $ is faithful, then  $\Phi$ is faithful.
  \end{enumerate}
\xprop 
\pf Put $\Gamma =\{d_i\}_{i \in K}$.  We have a natural application from $\Gamma$-sequences to $\Sigma$-sequences. This induces an application from $\Gamma$-fractions to $\Sigma$-fractions. This provides $\Phi$ by restriction. To prove (i), we observe that for each piece
of $\Sigma$-sequence of the kind $l_{d_i} \circ n$ with $i \in I \smallsetminus K$ and $n \in N_i$, by
assumption $n = d_i \circ f$ for some $f $ in $\calC$, so that $l_{d_i}$ disappears and the
remaining sequence is a $\Gamma$-sequence. To prove (ii), consider the commutative diagram of functors \[
\begin{tikzcd}
 \calC \left[ \left\{ (d_i)^{-1} \circ N_i \right\}_{i \in K} \right]  \ar[r] \ar[d, hook] & \calC \left[ \left\{ (d_i)^{-1} \circ N_i \right\}_{i \in I} \right]  \ar[d, hook] \\
\calC [ \Gamma^{-1} ] \ar[r, hook] & \calC [  \Sigma ^{-1} ] 
\end{tikzcd}.\]
The vertical arrows are faithful by Fact \ref{dildansloc}. The lower arrow is faithful by assumption. This shows that the upper arrow is faithful and finishes (ii).
\xpf 

\prop \label{indu} Let $\{ [N_j , d_j ]\} _{j \in J}$ be another center in $\calC$. Put $I' = I \sqcup J$. 
Proposition \ref{propdil} provides functors 
\[\Theta : \calC \to \calC \left[ \{ (d_i)^{-1} \circ N_i \} _{i \in I }  \right]\]  and 
\[\Theta ' : \calC \to \calC \left[ \{ (d_i)^{-1} \circ N_i \} _{i \in I' }  \right]. \]We observe that $\left\{ \left[ S_{\Theta (N_j) }^{\calC [ \{ (d_i)^{-1} \circ N_i \} _{i \in I }  ]} , \Theta (d_j)\right] \right\}_{j \in J} $ is a center in $\calC [ \{ (d_i)^{-1} \circ N_i \} _{i \in I }  ]$, so that we get a dilatation functor  
 \[ \beta : \calC [ \{ (d_i)^{-1} \circ N_i \} _{i \in I }  ] 
 \to 
 \calC [ \{ (d_i)^{-1} \circ N_i \} _{i \in I }  ] \left[ \left\{ (\Theta (d_j))^{-1} \circ S_{\Theta (N_j) }^{\calC [ \{ (d_i)^{-1} \circ N_i \} _{i \in I }  ]} \right\} _{j \in J }   \right] .\]
We observe that Proposition \ref{fi} provides a canonical functor
 \[\Phi :  \calC [ \{ (d_i)^{-1} \circ N_i \} _{i \in I }  ]\to \calC [ \{ (d_i)^{-1} \circ N_i \} _{i \in I' }  ] .\]
\begin{enumerate}
\item The functor $\Phi $ belongs to $Cat_{\calC [ \{ (d_i)^{-1} \circ N_i \} _{i \in I }  ] }^{\{ \Theta (d_j)\}_{j \in J }\text{-reg}}$.
\item The functor $\beta\circ \Theta$ belongs to $Cat_{\calC} ^{\{ d_i\}_{i \in I'}\text{-reg}}.$
\item \label{iiift} There exists a unique functor \[ \alpha : \calC [ \{ (d_i)^{-1} \circ N_i \} _{i \in I' }  ] \to   \calC [ \{ (d_i)^{-1} \circ N_i \} _{i \in I }  ] \left[ \left\{ (\Theta (d_j))^{-1} \circ S_{\Theta (N_j) }^{\calC [ \{ (d_i)^{-1} \circ N_i \} _{i \in I }  ]} \right\} _{j \in J }   \right] \] such that $\beta \circ \Theta  = \alpha \circ  \Theta' $.

\item \label{ivft} There exists a unique functor  \[ \alpha' :   \calC [ \{ (d_i)^{-1} \circ N_i \} _{i \in I }  ] \left[ \left\{ (\Theta (d_j))^{-1} \circ S_{\Theta (N_j) }^{\calC [ \{ (d_i)^{-1} \circ N_i \} _{i \in I }  ]} \right\} _{j \in J }   \right] \to  \calC [ \{ (d_i)^{-1} \circ N_i \} _{i \in I' }  ] \] such that $\Phi  = \alpha ' \circ \beta  $.

\item \label{vft}The functors $\alpha \circ \alpha ' $ and $\alpha ' \circ \alpha $ are identity functors.

\item \label{vift}There exists an isomorphism of categories 
\[    \calC [ \{ (d_i)^{-1} \circ N_i \} _{i \in I }  ] \left[ \left\{ (\Theta (d_j))^{-1} \circ S_{\Theta (N_j) }^{\calC [ \{ (d_i)^{-1} \circ N_i \} _{i \in I }  ]} \right\} _{j \in J }   \right] \isomto \calC [ \{ (d_i)^{-1} \circ N_i \} _{i \in I' }  ] .\] 
\end{enumerate}
\xprop 
\pf
\begin{enumerate}
\item  We have to prove that $\calC [ \{ (d_i)^{-1} \circ N_i \} _{i \in I' }  ] \to \calC [ \{ (d_i)^{-1} \circ N_i \} _{i \in I' }  ] [(\{\Phi ( \Theta (d_j)) \}_{j \in J })^{-1}]  $ is faithful. It is enough to observe that $ \Phi \circ \Theta = \Phi' $ and apply Fact \ref{obecjtreg} and Fact \ref{factr}.
\item Applying Fact \ref{dildansloc} twice, we get a faithful functor 
\begin{equation} \label{o1}\calC [ \{ (d_i)^{-1} \circ N_i \} _{i \in I }  ] [ \{ (\Theta (d_j))^{-1} \circ S_{\Theta (N_j) }^{\calC [ \{ (d_i)^{-1} \circ N_i \} _{i \in I }  ]} \} _{j \in J }   ] \to \calC [(\{d_i\}_{i \in I})^{-1}] [(\{[d_j]\}_{j \in J})^{-1}].\end{equation} Now we observe that \begin{equation} \label{o2} \calC [(\{d_i\}_{i \in I})^{-1}] [(\{[d_j]\}_{j \in J})^{-1}]=\calC [(\{d_i\}_{i \in I'})^{-1}]  \end{equation}  and that \begin{equation}\label{o3}\calC [ \{ (d_i)^{-1} \circ N_i \} _{i \in I }  ] [ \{ (\Theta (d_j))^{-1} \circ S_{\Theta (N_j) }^{\calC [ \{ (d_i)^{-1} \circ N_i \} _{i \in I }  ]} \} _{j \in J }   ] [ (\{ (\beta \circ \Theta ) (d_i) \}_{i \in I'})^{-1} ] = \calC [(\{d_i\}_{i \in I'})^{-1}] . \end{equation} 

Now (\ref{o1}), (\ref{o2}) and (\ref{o3}) together implies that the functor $\beta\circ \Theta$ belongs to $Cat_{\calC} ^{\{ d_i\}_{i \in I'}\text{-reg}}.$ 
\item By Theorem \ref{uni}, it is enough to prove that, for all $k $ in $I'$:
{\small \begin{equation}\label{eqazd} S_{(\beta \circ \Theta ) ( N_k)}^{\calC [ \{ (d_i)^{-1} \circ N_i \} _{i \in I }  ] [ \{ (\Theta (d_j))^{-1} \circ S_{\Theta (N_j) }^{\calC [ \{ (d_i)^{-1} \circ N_i \} _{i \in I }  ]} \} _{j \in J }   ]} \subset S_{(\beta \circ \Theta ) (d_k)}^{\calC [ \{ (d_i)^{-1} \circ N_i \} _{i \in I }  ] [ \{ (\Theta (d_j))^{-1} \circ S_{\Theta (N_j) }^{\calC [ \{ (d_i)^{-1} \circ N_i \} _{i \in I }  ]} \} _{j \in J }   ]} .\end{equation}}
Let us first take $k $ in $I$. By Proposition \ref{exc}, we have
\[
S_{\Theta (N_k) }^{\calC [ \{ (d_i)^{-1} \circ N_i \} _{i \in I }  ]} \subset S_{\Theta (d_k) }^{\calC [ \{ (d_i)^{-1} \circ N_i \} _{i \in I }  ]}
\]
therefore we get 
{\small \begin{equation} \label{eqzdf} S_{\beta \left(S_{\Theta (N_k) }^{\calC [ \{ (d_i)^{-1} \circ N_i \} _{i \in I }  ]}\right)}^{\calC [ \{ (d_i)^{-1} \circ N_i \} _{i \in I }  ] [ \{ (\Theta (d_j))^{-1} \circ S_{\Theta (N_j) }^{\calC [ \{ (d_i)^{-1} \circ N_i \} _{i \in I }  ]} \} _{j \in J }   ]} \subset S_{\beta \left(S_{\Theta (d_k) }^{\calC [ \{ (d_i)^{-1} \circ N_i \} _{i \in I }  ]}\right)}^{\calC [ \{ (d_i)^{-1} \circ N_i \} _{i \in I }  ] [ \{ (\Theta (d_j))^{-1} \circ S_{\Theta (N_j) }^{\calC [ \{ (d_i)^{-1} \circ N_i \} _{i \in I }  ]} \} _{j \in J }   ]} .\end{equation}} 

Now (\ref{eqzdf}) and Fact \ref{factbd} implies (\ref{eqazd}) for all $k \in I$.
Now let $k \in J $.
By Proposition \ref{exc}, we have 
{\small \begin{equation} \label{eqzdf3} S_{\beta \left(S_{\Theta (N_k) }^{\calC [ \{ (d_i)^{-1} \circ N_i \} _{i \in I }  ]}\right)}^{\calC [ \{ (d_i)^{-1} \circ N_i \} _{i \in I }  ] [ \{ (\Theta (d_j))^{-1} \circ S_{\Theta (N_j) }^{\calC [ \{ (d_i)^{-1} \circ N_i \} _{i \in I }  ]} \} _{j \in J }   ]} \subset S_{\beta \left(\Theta (d_k) \right)}^{\calC [ \{ (d_i)^{-1} \circ N_i \} _{i \in I }  ] [ \{ (\Theta (d_j))^{-1} \circ S_{\Theta (N_j) }^{\calC [ \{ (d_i)^{-1} \circ N_i \} _{i \in I }  ]} \} _{j \in J }   ]} .\end{equation}} 
Now (\ref{eqzdf3}) and Fact \ref{factbd} implies (\ref{eqazd}) for all $k \in J$.
Finally, (\ref{eqazd}) holds for all $k \in I'.$
\item By Theorem \ref{uni}, it is enough to prove that, for all $k$ in $J$, 
\begin{equation} \label{ejk}
S^{\calC [ \{ (d_i)^{-1} \circ N_i \} _{i \in I' }  ]}_{\Phi ( S^{\calC [ \{ (d_i)^{-1} \circ N_i \} _{i \in I }  ]}_{\Theta ( N_k)}) } \subset S^{\calC [ \{ (d_i)^{-1} \circ N_i \} _{i \in I' }  ]} _{\Phi ( \Theta (d_k) )} .
\end{equation}
Using that $\Theta ' = \Phi \circ \Theta$,  Proposition \ref{exc} shows that, for all $k$ in $J$,  
\begin{equation} \label{ejkp}
S^{\calC [ \{ (d_i)^{-1} \circ N_i \} _{i \in I' }  ]}_{\Phi (\Theta ( N_k) )} \subset S^{\calC [ \{ (d_i)^{-1} \circ N_i \} _{i \in I' }  ]} _{\Phi ( \Theta (d_k) )} .
\end{equation}
Now (\ref{ejkp}) and Fact \ref{factbd} implies (\ref{ejk}) for all $k \in J$, as required.
\item This is an immediate consequence of Theorem \ref{repre}.
\item This is an immediate consequence of (\ref{vft}).
\end{enumerate} 
\xpf 

\rema
One can prove Proposition \ref{indu} (\ref{vift}) directly using explicit computations on fractions. 
However, Proposition \ref{indu} (\ref{iiift}) (\ref{ivft}) and (\ref{vft}) provides stronger uniqueness assertions.
\xrema

\fact \label{555}
Let $H: \calA \to \calB $ be a functor. Let $N,N'$ be two sieves of $A$. Then \[S^{\calB} _{H (N) } \cup S^{\calB} _{H (N')} = S_{H (N \cup N')} ^{\calB} .\]
\xfact 
\pf
The inclusion $\subset $ is immediate. Let $x $ be in $S_{H (N \cup N')} ^{\calB}$, then $x= H(n) \circ y $ with $n \in N \cup N'$. So $x$ belongs to $S^{\calB} _{H (N) } \cup S^{\calB} _{H (N')}$.
\xpf 

\prop For all $i \in I$, let $N_i '$ be another sieve over $cod(d_i)$ and let $N_i''$ be the union of $N_i$ and $N_i'$.
Then $\calC [ \{ (d_i)^{-1} \circ N_i \} _{i \in I }, \{ (d_i)^{-1} \circ N_i' \} _{i \in I }  ] $ identifies with $\calC [ \{ (d_i)^{-1} \circ N_i'' \} _{i \in I } ]$ (more precisely there are unique $\calC$-functors from each category to the other and these morphisms are isomorphisms).
\xprop 
\pf
We observe that $\Theta _1: \calC \to \calC [ \{ (d_i)^{-1} \circ N_i \} _{i \in I }, \{ (d_i)^{-1} \circ N_i' \} _{i \in I }  ] $ and  $\Theta _2 : \calC \to \calC [ \{ (d_i)^{-1} \circ N_i'' \} _{i \in I } ]$ belongs to $Cat_{\calC}^{\Sigma\text{-reg}}$ by Fact \ref{obecjtreg}. 
For $k \in I$, by  Proposition \ref{exc} we have \[S^{\calC [ \{ (d_i)^{-1} \circ N_i \} _{i \in I }, \{ (d_i)^{-1} \circ N_i' \} _{i \in I }  ]}_{\Theta _1 (N_k)}  \cup S^{\calC [ \{ (d_i)^{-1} \circ N_i \} _{i \in I }, \{ (d_i)^{-1} \circ N_i' \} _{i \in I }  ]}_{\Theta _1 (N_k')}  \subset S^{\calC [ \{ (d_i)^{-1} \circ N_i \} _{i \in I }, \{ (d_i)^{-1} \circ N_i' \} _{i \in I }  ]}_{\Theta _1 (d_k)} . \]
So by Fact \ref{555} we have 
\[S^{\calC [ \{ (d_i)^{-1} \circ N_i \} _{i \in I }, \{ (d_i)^{-1} \circ N_i' \} _{i \in I }  ]}_{\Theta _1 (N_k \cup N_k') }   \subset S^{\calC [ \{ (d_i)^{-1} \circ N_i \} _{i \in I }, \{ (d_i)^{-1} \circ N_i' \} _{i \in I }  ]}_{\Theta _1 (d_k)} . \]
So by Theorem \ref{uni}, we get a unique $\calC$-functor 
\[ \alpha : \calC [ \{ (d_i)^{-1} \circ N_i'' \} _{i \in I } ] \to  \calC [ \{ (d_i)^{-1} \circ N_i \} _{i \in I }, \{ (d_i)^{-1} \circ N_i' \} _{i \in I }  ]. \]

For $k\in I$, by Proposition \ref{exc}, we have 
\[ S^{\calC [ \{ (d_i)^{-1} \circ N_i'' \} _{i \in I }\} _{i \in I }  ]}_{\Theta _2 (N_k'')}   \subset S^{\calC [ \{ (d_i)^{-1} \circ N_i \} _{i \in I }]}_{\Theta _2 (d_k)} . \]

Since $N_k,N_k' \subset N_k'', $ by Theorem \ref{uni}, we get a unique $\calC$-functor 
\[ \alpha' :  \calC [ \{ (d_i)^{-1} \circ N_i \} _{i \in I }, \{ (d_i)^{-1} \circ N_i' \} _{i \in I }  ] \to \calC [ \{ (d_i)^{-1} \circ N_i'' \} _{i \in I } ] . \]
Proposition \ref{repre} implies that $\alpha \circ \alpha '$ and $\alpha ' \circ \alpha$ are identity morphisms. 
\xpf

\section{Codilatations of categories} \label{codi} \label{leftright}
We used fractions of morphisms $(d_i)^{-1} \circ N_i$, the construction given in this note also makes sense for fractions $V_i \circ  (d_i)^{-1}$ (each $V_i$ is now a cosieve), pictorially:
\[\begin{tikzcd}
X \ar[dr, "d_i"] \ar[rr, "v" ] && Y\\
& Z  \ar[ru, dotted, "\exists !"]&
\end{tikzcd} \forall i \in I  , \forall v \in V_i.\]

In this section, we introduce this construction formally. We define codilatations using dilatations and opposite categories. Let $\calC $ be a category.

\defi
 A cocenter $\{[V_i , d_i]\}_{i \in I}$ in $\calC$ is a collection of pairs $[V_i , d_i]$ where $d_i$ is a morphism of $\calC$ and $V_i$ is a cosieve from $\dom (d_i)$.
\xdefi 
We now fix a cocenter $\{[V_i , d_i]\}_{i \in I}$.
\fact
In $\calC ^{op}$, for all $i $ in $I$, $V_i$ can be regarded as sieve over $cod^{\calC ^{op}}(d_i)= dom ^{\calC}(d_i)$. In particular $\{[V_i , d_i]\}_{i \in I}$ can be regarded as a center in $\calC ^{op}$.
\xfact

\pf
By definition, the collection of morphisms of a cosieve is a sieve in the opposite category. 
\xpf

\defi \begin{sloppypar}
We put $\calC \left[ \{V_i \circ (d_i)^{-1} \} _{i \in I }\right] = \left( \calC ^{op} \left[ \{ (d_i)^{-1} \circ V_i  \}_{i \in I } \right] \right) ^{op} $. The category $\calC \left[ \{V_i \circ (d_i)^{-1} \} _{i \in I }\right] $ is called the codilatation of $\calC$ with cocenter $\{[V_i , d_i]\}_{i \in I}$.\end{sloppypar}
\xdefi 

\fact \label{sievcosievop}
Let $\calA$ be a category. Let $E$ be a collection of morphisms of $\calA$. Then $CoS^{A}_E = S^{A^{op}}_{E}$.
\xfact
 \pf
 $ CoS^{\calA}_{E} = \{ x \circ _{\calA} e | dom ^\calA (x) = cod ^\calA (e) \} = \{ e \circ _{\calA ^{op}} x | dom ^{\calA ^{op}} (e) = cod ^{\calA}(x)\} = S_{E}^{\calA^{op}}$.
 \xpf 
\prop 
\begin{enumerate}
 \item We have a canonical functor $\Upsilon : \calC \to \calC \left[ \{V_i \circ (d_i)^{-1} \} _{i \in I }\right]$.
 \item We have a canonical faithful functor $\calC \left[ \{V_i \circ (d_i)^{-1} \} _{i \in I }\right] \to \calC \left[ (\{d_i \} _{i \in I })^{-1}\right].$
 \item The functor $\Upsilon $ is an object in the category $Cat _{\calC}^{\{d_i \} _{i \in I } \text{-reg}}$.
 \item The functor $\Upsilon $ represents the covariant functor $Cat _{\calC}^{\{d_i \} _{i \in I } \text{-reg}}$ to $Set$ 
given by
\begin{equation}
(\calC \xrightarrow{F} \calD) \;\longmapsto\; \begin{cases}\{*\}, \; \text{ if $CoS^{\calD}_{F(v_i)} \subset CoS^{\calD}_{F(d_i)}$ for all $i $ in $I$;}\\ \varnothing,\;\text{else.}\end{cases}
\end{equation}

\end{enumerate}
\xprop 

\pf \begin{sloppypar}
\begin{enumerate}
\item The canonical functor $\calC ^{op} \to  \calC ^{op} \left[ \{ (d_i)^{-1} \circ V_i  \}_{i \in I } \right] $ of Proposition \ref{propdil} induces a canonical functor $ \Upsilon :\calC \to  \left( \calC ^{op} \left[ \{ (d_i)^{-1} \circ V_i  \}_{i \in I } \right] \right) ^{op}$. 
\item Fact \ref{dildansloc} provides a faithful functor $\calC ^{op} \left[ \{ (d_i)^{-1} \circ V_i  \}_{i \in I } \right] \to \calC ^{op} [(\{d_i\}_{i \in i})^{-1}] $.
 Since $\calC ^{op} [(\{d_i\}_{i \in i})^{-1}] = \left (\calC  [(\{d_i\}_{i \in i})^{-1}]) \right) ^{op}$ (e.g. using the explicit descriptions with fractions),
 we get the desired faithful functor $\left( \calC ^{op} \left[ \{ (d_i)^{-1} \circ V_i  \}_{i \in I } \right] \right)^{op} \to \calC [(\{d_i\}_{i \in i})^{-1}]$.
 \item Fact \ref{obecjtreg} implies that $\calC ^{op} \to  \calC ^{op} \left[ \{ (d_i)^{-1} \circ V_i  \}_{i \in I } \right] $ is an object in the category $Cat _{\calC^{op}}^{\{d_i \} _{i \in I } \text{-reg}}$. This implies that the functor $\Upsilon $ is an object in the category $Cat _{\calC}^{\{d_i \} _{i \in I } \text{-reg}}$.
 \item Combine Proposition \ref{repre} and Fact \ref{sievcosievop}.
\end{enumerate} \end{sloppypar}
\xpf 
\section{Some examples and remarks}
\subsubsection{Universal property of localizations of categories}
Let $\calC$ be a category and let $\Sigma $ be a class of morphisms of $\calC$. 
Fact \ref{locasdil} shows that the localization $\calC [\Sigma^{-1}] $ is equal to the dilatation $\calC [ \{ d ^{-1} \circ N_d \} _{d \in \Sigma} ]$ where $N_d = S_{Id_{cod(d)}}^{\calC}$. Obviously, $\calC[\Sigma^{-1}] \big[ [\Sigma]^{-1} \big]=\calC[\Sigma^{-1}] $ and so $\calC[\Sigma^{-1}]$ belongs to $Cat _{\calC}^{\Sigma \text{-reg}}$. Now if $F : \calC \to \calD$ is another $\calC$-category, such that $F(\Sigma)$ is made of isomorphisms, again $\calD= \calD [ F(\Sigma)^{-1}]$ so that $\calD$ also belongs to $Cat _{\calC}^{\Sigma \text{-reg}}$. Moreover in this case it is obvious that for any $d \in \Sigma $,  we have \[ S^{\calD}_{F(N_d)} \subset  S^{\calD}_{F(d)}.\] So by Theorem \ref{uni}, $F$ factors uniquely through $\calC'$. So the universal property of dilatations generalizes the universal property of localizations.

\subsubsection{Dilatations of commutative rings an semirings} \label{examplering}

\prop \label{Ma} Let $A$ be a commutative unital ring. Let $I$ be a set and let $\{ [ M_i, a_i] \} _{i \in I } $ be a center in $A$ in the sense of \cite{M}.  Let $A [ \{ \frac{M_i}{a_i} \} _{i \in I } ]$ be the dilatation of $A$ with center $\{ [ M_i, a_i] \} _{i \in I } $ as in \cite{M}. Let $\calC$ be the category attached to $A$. In particular $a_i$ is a morphism of $\calC$ and $M_i$ is a sieve for all $i$ in $I$. We have $Ob \calC = \{\bullet\}$, a singleton. Let $\calC [ \{ \frac{M_i}{a_i} \} _{i \in I } ]$ be the category attached to the ring $A[ \{ \frac{M_i}{a_i} \} _{i \in I } ]$.
  We have an identification of $\calC$-functors $\calC [\{ (a_i)^{-1} \circ M_i \}_{i \in I} ] =  \calC [ \{ \frac{M_i}{a_i} \} _{i \in I } ]$.
\xprop 
\pf  
Let $\Phi $ be the functor  $\calC [\{ (a_i)^{-1} \circ M_i \}_{i \in I} ] \to \calC [ \{ \frac{M_i}{a_i} \} _{i \in I } ]$ given by $[\bullet \xrightarrow{m} \bullet \xrightarrow{l_{a_i}} \bullet ] \mapsto \frac{m}{a_i} $. It is well-defined, surjective, injective and provides the desired identification.
\xpf 

This shows that dilatations of categories generalize dilatations of rings. 
As noticed in \cite{M}, dilatations of commutative semirings make sense. The same argument as before shows that dilatations of semirings also provide examples of dilatations of categories.

\subsubsection{Dilatations of monoids}

Let $\calC$ be a small category with a single object, i.e. a (not necessarily commutative) monoid $M$. Let $\{[N_i, d_i]\}_{i \in I }$ be a center of $\calC$. Then $\calC [\{(d_i)^{-1} \circ  N_i \}_{i \in I }]$ is a category with a single object endowed with a functor $\calC \to \calC [\{(d_i)^{-1} \circ N_i \}_{i \in I }]$. In other words, $\calC [\{(d_i)^{-1} \circ N_i \}_{i \in I }]$ is an $M$-monoid $M'$. We refer to this construction as dilatations of monoids. This generalizes localizations of monoids. Codilatations of monoids also make sense.

\subsubsection{Dilatations of pre-additive categories and non-commutative rings}

Let $\calC$ be a pre-additive category, in general a dilatation of $\calC$ is not pre-additive as we know that this fails already for localization, cf. e.g. \cite{A}. Nevertheless, as in the case of localizations, it should be possible to study dilatations of pre-additive categories once we have an adapted calculus of fractions. This is  related to investigate dilatations of non-commutative rings.

\subsubsection{A characterization of dilatations of rings that does not hold for categories}

In the context of dilatations of rings, \cite[Fact 2.14]{M} shows that for any sub-$A$-algebra $B$ of a localization $A[(\{a_i\}_{i \in I})^{-1}]$, there is a multi-center in $A$ such that $B $ identifies with the associated dilatation of $A$. Here we explain that this characterization does not extend to dilatations of categories. 

\fact \label{contrexemplecat}There exists a category $\calC $, a collection of morphisms $\Gamma$, a functor $ \calC \to \calD$  that is the identity on objects, and a faithful functor $\calD \to \calC [\Gamma ^{-1}]$ such that 
\[\begin{tikzcd} \calC \ar[rd]\ar[rr]&& \calC [\Gamma ^{-1}]  \\ &\calD \ar[ru]& \end{tikzcd}\] commutes and such that $\calC \to \calD$ is not isomorphic to $\calC \to \calC [\{(d_i)^{-1} \circ N_i\}_{i \in I} ]$ for all centers $\{[N_i,d_i]\}_{i\in I}$.
\xfact 

\pf
Let $\calC$ be the category with two objects $X $ and $Y$ and whose morphisms are as follows \[Hom _{\calC} (X,X)= Id _X, Hom _{\calC} (Y,Y) = Id_Y, Hom _{\calC} (Y,X)= \emptyset\text{ and }Hom _{\calC} (X,Y)= \{a,b\}\]

Let $\Gamma $ be the collection of morphisms of $\calC$ given by 
\[ \Gamma = \{  b \}.\]

Then $\calC [ \Gamma ^{-1} ]$ is the category whose objects are $X$ and $Y$ and whose morphisms are as follows
\begin{align*} Hom_{\calC [\Gamma ^{-1}]}(X,X)& = \{Id_X, b^{-1} a ,b^{-1} a b^{-1} a  , \ldots, (b^{-1} a )^n, \ldots  (n \in \bbN) \}, \\ Hom_{\calC [\Gamma ^{-1}]}(Y,Y) & = \{ Id_Y , ab^{-1}, ab^{-1} a b^{-1}, \ldots (ab^{-1})^, \ldots (n \in \bbN) \},  \\Hom_{\calC [\Gamma ^{-1}]}(X,Y)& = \{ b, a, ab^{-1}a, ab^{-1}a b^{-1} a , \ldots, a(b^{-1}a)^n, \ldots (n \in \bbN) \} ,\\
Hom_{\calC [\Gamma ^{-1}]}(Y,X) &= \{ b^{-1}, b^{-1}ab^{-1},  b^{-1} a b^{-1}ab^{-1}, \ldots, b^{-1} (ab^{-1})^n, \ldots (n \in \bbN) \} . \end{align*}

Now let $\calD$ be the subcategory of $\calC [\Gamma ^{-1}]$ whose objects are $X$ and $Y$ and whose morphisms are as follows 

\begin{align*} Hom_{\calD}(X,X)& = Id_X, \\ Hom_{\calD}(Y,Y)& = Id_Y , \\Hom_{\calD}(X,Y)& = \{ b, a,  ab^{-1}a  \} ,\\
Hom_{\calD}(Y,X) &= \varnothing . \end{align*}

Then we have canonical functors $\calC \to \calD $ and $\calD \to \calC [ \Gamma ^{-1}]$. We claim that there is no center $ \{[N_i , d_i ] \}_{i \in I }$ such that $\calD$ identifies with the associated dilatation.
To prove this claim, we chose a center  $ \{[N_i , d_i ] \}_{i \in I }$ in $\calC$, put $\Sigma = \{d_i\}_{i \in I}$ and exhaustively distinguish two cases:
\begin{enumerate}
\item if ($a \in \Sigma$ and ($Id_Y \in N_a $ or $b \in N_a$)) or ($b \in \Sigma $ and ($Id_Y \in N_b $ or $a \in N_b$)), then $\# Hom _{ \calC'} (X,Y) $ is infinite,
\item  if ($a \in \Sigma \Rightarrow N_a \subset \{ a\}$) and ($b \in \Sigma \Rightarrow N_b \subset \{ b \} $), then $\calC' = \calC$. \end{enumerate}
In all cases $ \calC' \ne  \calD $.
\xpf

\bigbreak\bigbreak
\noindent \par
\noindent Einstein Institute of Mathematics \\
Edmond J. Safra Campus\\
The Hebrew University of Jerusalem\\
Givat Ram. Jerusalem, 9190401, Israel \par
\noindent E-mail address: \texttt{arnaud.mayeux@mail.huji.ac.il}
\end{document}